\documentclass[12pt]{amsart}
\usepackage{amsfonts, cite, amssymb, setspace, amsmath, multicol, epsfig, amsthm, enumerate, hyperref, setspace}
\usepackage[utf8]{inputenc}
\usepackage[margin=1.2in]{geometry}
\makeatletter
\renewcommand{\@secnumfont}{\bfseries}
\makeatother

\newtheorem{theorem}{Theorem}[section]
\newtheorem*{thmV}{Theorem 1.6}

\theoremstyle{definition}

\newtheorem{example}[theorem]{Example}

\newtheorem{corollary}[theorem]{Corollary}
\newtheorem*{acknowledgement}{Acknowledgement}

\theoremstyle{remark}
\newtheorem*{remark}{Remark}

\numberwithin{equation}{section}
\usepackage{ amssymb }

\begin{document}

\title{Moonshine Modules and a question of Griess}

\author[V.\,M. Aricheta]{Victor Manuel Aricheta}
\address{Department of Mathematics and Computer Science, Emory University}
\email{victor.manuel.aricheta@emory.edu}

\author[L. Beneish]{Lea Beneish}
\address{Department of Mathematics and Computer Science, Emory University}
\email{lea.beneish@emory.edu}


\begin{abstract} We consider the situation in which a finite group acts on an infinite-dimensional graded module in such a way that the graded trace functions are weakly holomorphic modular forms. Under a mild hypothesis we completely describe the asymptotic module structure of the homogeneous subspaces. As a consequence we find that moonshine for a group gives rise to partial orderings on its irreducible representations. This serves as a first answer to a question posed by Griess. In particular, we show that our hypothesis holds for umbral moonshine and for automorphism groups of certain vertex operator algebras.
\end{abstract}

\maketitle



\section{Introduction}

In mathematics, the term moonshine is used to reference a situation in which distinguished modular or mock modular forms serve as graded trace functions for the action of a finite group on some infinite-dimensional module. In the late 1970s, Conway and Norton coined the phrase \emph{monstrous moonshine} to describe the (then-conjectural) relationship between the Monster group $\mathbb{M}$ and modular functions\cite{MR554399}. More precisely, monstrous moonshine asserts the existence of a graded infinite-dimensional $\mathbb{M}$-module
	\[ V^{\natural} = V^{\natural}_{-1} \oplus V^{\natural}_{1} \oplus V^{\natural}_{2} \oplus \cdots \]
whose graded trace functions are generators of function fields of genus zero quotients of the complex upper half-plane. The \emph{monstrous moonshine module} $V^{\natural}$ was constructed by Frenkel, Lepowsky, and Meurman in 1988  \cite{MR781381}. The fact that the graded trace functions are generators of function fields was established by Borcherds in 1992 via the theory of vertex operator algebras and generalized Kac-Moody algebras \cite{MR1172696}.  

More recently in 2010, in their study of the elliptic genus of a K3 surface, string theorists Eguchi, Ooguri, and Tachikawa discovered a moonshine phenomenon relating the largest Mathieu group $M_{24}$ with a certain mock modular form\cite{EH}. Analogous to monstrous moonshine, we thus have \emph{Mathieu moonshine}, which posits the existence of a graded infinite-dimensional $M_{24}$-module  $K^\natural = \bigoplus K^{\natural}_n$, called the \emph{Mathieu moonshine module}, whose graded trace functions are mock modular forms. The existence of $K^\natural$ was proven by Gannon \cite{Gannon}. We now know that Mathieu moonshine is but one of the 23 cases of moonshine referred to collectively as \emph{umbral moonshine}, which was discovered by Cheng, Duncan, and Harvey \cite{CDH} and proven later by Duncan, Griffin, and Ono \cite{proof}. 



Duncan, Griffin, and Ono revisited monstrous moonshine and explored the structure of the homogeneous subspaces $V^{\natural}_{n}$ of the monstrous moonshine module \cite{MR3375653}. They found that as $n \to \infty$, the subspaces $V^{\natural}_{n}$ tend to a multiple of the regular representation of $\mathbb{M}$. To explain this, let $M^{(\mathbb{M})}_1, \ldots, M^{(\mathbb{M})}_{194}$ be the irreducible representations of $\mathbb{M}$, ordered as in the ATLAS \cite{conway1985atlas}. Write
	\[ V^{\natural}_{n} = \text{m}_1(n)M^{(\mathbb{M})}_1 \oplus \cdots \oplus \text{m}_{194}(n)M^{(\mathbb{M})}_{194}. \]
For $i = 1$, $\ldots$, $194$, they showed that 
	\begin{equation}\label{monstermi} \text{m}_i(n) \sim \dfrac{e^{4 \pi \sqrt{n}}}{\sqrt{2}|\mathbb{M}|n^{3/4}} \, \dim M^{(\mathbb{M})}_i \end{equation}
as $n \to \infty$. To show this, they derived an exact series formula for $\text{m}_i(n)$, and the asymptotic above is obtained by isolating the dominant term of the series. In view of this result, Griess posed the following question (cf. Problem 10.9. in \cite{MR3375653}):

\begin{quote}
\textbf{Griess' Question.} \emph{If we write each homogeneous subspace of each moonshine module, particularly the moonshine module $V^{\natural}$, as the sum of a free part (free over the group ring of $\mathbb{M}$) and a non-free part, then the non-free part tends to 0 (relative to the free part) as $n\to \infty$. Is there something to be learnt from an analysis of the non-free parts?}
\end{quote}

A step in this direction is given by Larson, who found asymptotic formulas for the non-free parts of $V_n^{\natural}$ \cite{Hannah}. We point out here that by closely reading the asymptotic formulas, one sees that the non-free parts of $V_n^{\natural}$ tend to a representation of $\mathbb{M}$ whose irreducible components do not include $M^{(\mathbb{M})}_{16}$ and $M^{(\mathbb{M})}_{17}$. We will show that this extends even further (cf. Theorem 1.2). To explain this, we first introduce a definition.

Let $(K_n)$ be a sequence of finite-dimensional representations of a finite group $G$, and suppose $c_g(n):= \text{tr}(g|K_n) \in \mathbb{R}$ for all $g\in G$ and all $n$. We say that the sequence $(K_n)$ has \emph{dominant identity trace} if for every $g  \in G$ that is not equal to the identity element $e$, we have $ c_g(n) = o(c_e(n))$ as $n\rightarrow \infty$. Examples of such sequences are the sequences $(V^{\natural}_n)$ and $(K^{\natural}_n)$. Certain vertex operator algebras, for which the monstrous moonshine module is a special case, also have sequences of homogeneous subspaces that have dominant identity trace (cf. Theorem 1.6).

The following theorem shows that if a sequence $(K_n)$ has dominant identity trace, then the subspaces $K_n$ tend to a multiple of the regular representation of $G$. Thus, Griess' question makes sense for these sequences.

\begin{theorem}\label{thm1}
Let $e$ be the identity element in $G$, and let $M_1, \ldots, M_s$ be the irreducible representations of $G$. Write
	\[K_{n}^{} = \emph{m}_{1}(n)M_1 \oplus \emph{m}_{2}(n)M_2 \oplus \cdots \oplus \emph{m}_{s}(n)M_s.\]
If $(K_n)$ has dominant identity trace, then
	\[ \emph{m}_{i}(n) \sim \dfrac{1}{|G|} \dim K_n \dim M_i \]  as $n \to \infty$. Therefore, \[\lim_{n\to \infty} \dfrac{\emph{m}_i(n)}{\sum_{j=1}^{s}\emph{m}_j(n)}=\dfrac{\dim M_i}{\sum_{j=1}^s \dim M_j} .\] 
\end{theorem}

\vspace{0.1in}

\begin{remark} In some cases, $\dim K_n$ has known asymptotics in terms of simple functions.  For example, in the case of monstrous moonshine the asymptotics for the coefficients of the $j$-function yield \eqref{monstermi}. As another example, if $(K^{\natural}_n)$ is the sequence of homogeneous subspaces of the Mathieu moonshine module, then $\dim K_n^{\natural}$ may be written as a Rademacher series (cf. Section 3.1). This yields \[ \dim K^{\natural}_n \sim \frac{4}{\sqrt{8n-1}}\, \text{exp}\left(\frac{\pi\sqrt{8n-1}}{2}\right)\]
as $n\to \infty$. Therefore if $M^{(M_{24})}_i$ denotes an irreducible representation of $M_{24}$ and $\text{m}_i(n)$ denotes the multiplicity of $M^{(M_{24})}_i$ in $K^{\natural}_n$, then
\[ \text{m}_i(n) \sim \dfrac{4e^{\pi\sqrt{8n-1}/2}}{|M_{24}|\sqrt{8n-1}} \, \dim M^{(M_{24})}_i \]
as $n \to \infty$.
\end{remark}

\vspace{0.1in}


The statement of Theorem \ref{thm1}, that $K_{n}$ tends to a multiple of the regular representation of $G$, is obtained by an analysis which involves only the identity element of $G$. By performing an analysis which includes all the other elements of $G$, we find other representations of $G$---sensitive to some initial condition---which we can view as natural analogues of the regular representation. More precisely, we have the following result.

\begin{theorem}\label{thm2} Let $(K_n)$ be a sequence of finite-dimensional representations of a finite group $G$ and suppose $c_g(n):= \text{tr}(g|K_n) \in \mathbb{R}$ for all $g\in G$ and all $n$. Let $(n_i)$ be a sequence of integers such that given $g\in G$, the signs $\text{sgn}(c_g(n_i))$ are independent of $i$. If $(K_n)$ has dominant identity trace, then there exist $G$-modules $L_1, L_2, \ldots, L_t$ (depending on the signs $\text{sgn}(c_g(n_i)))$ where
\begin{itemize}
\item $L_1$ is the regular representation of $G$, and
\item the irreducible components of $L_{i+1}$ form a subset of the irreducible components of $L_i$ (for $1\leq i < t$), 
\end{itemize}
such that
\begin{enumerate}
\item for some nonnegative integer-valued functions $r_1(n_i), \ldots, r_t(n_i)$ and $G$-module $L_{\epsilon}(n_i)$ with bounded multiplicity functions, we have the decomposition \[ K_{n_i} = r_1(n_i)L_1 \oplus r_2(n_i)L_2 \oplus \cdots \oplus r_t(n_i)L_t \oplus L_{\epsilon}(n_i), \]
\item and the module $K_{n_i} \oplus (-r_1(n_i))L_1 \oplus \cdots \oplus (-r_l(n_i))L_l$ tends to a multiple of the representation $L_{l+1}$ (for $1 \leq l \leq t-1$) as $i \to \infty$.
\end{enumerate}
 \end{theorem}

This result shows that the representations $L_j$'s have a curious property that they are expressed in terms of fewer and fewer irreducible representations of $G$. In other words, by looking at the sequence $L_1, L_2, \ldots$, we find that the irreducible representations disappear in some order. Thus we find that moonshine for a group naturally equips its irreducible representations with partial orders. This is the part that speaks to Griess' question.

\begin{example} Let $H^{(M_{24})}_{g}(\tau)$ be the Mathieu moonshine graded trace functions (cf. Section 3.1). Consider the action of $A_5 \subset M_{24}$ on the Mathieu moonshine module $K^\natural = \bigoplus K_n^\natural$ with the following graded trace functions:
\begin{align*}
H_{1A}^{(A_5)} &= H_{1A}^{(M_{24})}\\
H_{2A}^{(A_5)} &= H_{2A}^{(M_{24})}\\
H_{3A}^{(A_5)} &= H_{3A}^{(M_{24})}\\
H_{5A}^{(A_5)} &= H_{5A}^{(M_{24})}\\
H_{5B}^{(A_5)} &= H_{5A}^{(M_{24})}.
\end{align*}
Let $M^{(A_5)}_1, \ldots, M^{(A_5)}_5$ be the irreducible representations of $A_5$, labelled as in GAP's SmallGroup library \cite{GAP4}. If $n_i=10+30i$, then as $i \to \infty$ the discussion in Section 2 shows that $K^{\natural}_{n_i}$ naturally decomposes as an $A_5$-module into
\[ K^{\natural}_{n_i} = r_1(n_i)L_1 \oplus r_2(n_i)L_2 \oplus r_3(n_i)L_3 \oplus L_{\epsilon}(n_i)\]
where: $L_1$ is the regular representation of $A_5$; $L_2$ is a representation of $A_5$ whose irreducible decomposition is in terms of $M^{(A_5)}_1$, $M^{(A_5)}_4$, and $M^{(A_5)}_5$; $L_3$ is a representation of $A_5$ whose irreducible decomposition is in terms of $M^{(A_5)}_1$ and $M^{(A_5)}_5$; and $L_{\epsilon}$ have bounded multiplicity functions. Hence, moonshine on $A_5$ gives us the partial ordering of (blocks of) irreducible modules: $\{M^{(A_5)}_2,M^{(A_5)}_3\}$, followed by $\{M^{(A_5)}_4\}$, and then by $\{M^{(A_5)}_1,M^{(A_5)}_5\}$. As another example, for $n_i= 21 +{30}i$, moonshine on $A_5$ gives us the partial ordering of (blocks of) irreducible modules: $\{M^{(A_5)}_1\}$, followed by $\{M^{(A_5)}_5\}$, and then by $\{M^{(A_5)}_2,M^{(A_5)}_3,M^{(A_5)}_4\}$.

An example where the ordering of the irreducible modules does not depend on the congruence class of $n$ is the following: Consider the action of $S_3\subset \mathbb{M}$ on the monstrous moonshine module $V^\natural$ with graded trace functions
\begin{align*}
H_{1A}^{(S_3)} &= T_{1A}\\
H_{2A}^{(S_3)} &= T_{2A}\\
H_{3A}^{(S_3)} &= T_{3A}
\end{align*}
where $T_{g}$ is the hauptmodul corresponding to $g\in \mathbb{M}$ via monstrous moonshine. This gives the ordering $\{ M^{(S_3)}_2 \}$, $ \{M^{(S_3)}_3 \}$, $\{M^{(S_3)}_1\}$, where $M^{(S_3)}_1$, $M^{(S_3)}_2$, $M^{(S_3)}_3$ are the irreducible representations of $S_3$ ordered by increasing dimension (starting with the trivial representation).
\end{example}

Griess' question is open ended, and the quantitative answer we offer is one of many partial answers to this problem. In fact, when Griess posed the problem, the original intention was to know whether a complete understanding of the asymptotics could provide clues to a richer algebraic structure surrounding the group and the graded module for it.

In the course of proving Theorem \ref{thm2}, we have also obtained the asymptotics for the multiplicities of the irreducible components of the non-free parts of $K_n$. We record this in the following theorem. Here, the set $\mathcal{C}_2$ is the collection of conjugacy classes of $G$ whose corresponding $c_g(n)$'s have the second fastest growth. (See the Proof of Theorems 1.2 and 1.4, in particular (\ref{C2}), for a precise definition of $\mathcal{C}_2$.)

\begin{theorem}\label{thm3}
Let $M_1, \ldots, M_s$ be the irreducible representations of a finite group $G$ and let $\chi_1, \ldots, \chi_s$ be their respective characters. Let $(K_n)$ be a sequence of $G$-modules that has dominant identity trace. Suppose $c_g(n):= \text{tr}(g|K_n) \in \mathbb{R}$ for all $g\in G$ and all $n$. Denote by $K_{n}'$ the non-free part of $K_{n}$ and write
	\[ K_{n}' = \emph{m}_1'(n)M_1 \oplus \emph{m}_2'(n)M_2 \oplus \cdots \oplus \emph{m}_s'(n)M_s.\]
Suppose that $(n_j)$ is a sequence of integers such that given $g \in \mathcal{C}_2$, the signs $\text{sgn}(c_g(n_j))$ are independent of $j$. Then as $j \to \infty$
	\[ \emph{m}_i'(n_j) \sim \dfrac{1}{|G|} \sum_{[g] \in \mathcal{C}_2} |[g]|f_i'(g)c_{g}(n_j). \]
Here
	\[ f_i'(g) := \overline{\chi_i(g)} -  \dfrac{\dim M_i}{\dim M_{j'}}\overline{\chi_{j'}(g)} \]
where $j'$ is a $j$ that minimizes
	\[ \dfrac{\sum\limits_{g \in \mathcal{C}_2} |[g]|\overline{\chi_j(g)}\emph{sgn}(c_g(n))}{\dim M_j} \]
as $n \to \infty$.
\end{theorem}

We discuss in Section $3$ some cases where our results apply. First, we show that our results apply to umbral moonshine. In umbral moonshine the trace functions are mock modular forms which have coefficients that can be expressed in terms of Kloosterman sums weighted by Bessel functions. Thus the coefficients of the trace functions have known asymptotics for which our theory applies. As an explicit example, we use Theorem $1.4$ to obtain the asymptotics of the multiplicities of the non-free part of the Mathieu moonshine module, which we record here as a corollary.

\begin{corollary}\label{cor5} Let $K^{\natural} = \bigoplus K^{\natural}_n$ be the Mathieu moonshine module, and let $K_{n}'$ be the non-free part of $K^{\natural}_{n}$. Let $M_1^{(M_{24})},\ldots, M_{26}^{(M_{24})}$ be the irreducible representations of $M_{24}$, and let $\chi_1, \ldots, \chi_{26}$ be their respective characters. Write
	\[ K_{n}' = \text{m}_1'(n)M^{(M_{24})}_1 \oplus \text{m}_2'(n)M^{(M_{24})}_2 \oplus \cdots \oplus \text{m}_{26}'(n)M^{(M_{24})}_{26}.\]
Then
\begin{align*}
	\text{m}_i'(n) &\sim (-1)^{n+1}\frac{\sqrt{2} \, e^{\frac{\pi}{4}\sqrt{8n-1}}}{\sqrt{8n-1}}\left(\frac{|2A|  }{|M_{24}|}\left( \chi_i(2A) - \frac{\chi_j(2A)}{\dim M^{(M_{24})}_j}\dim  M^{(M_{24})}_i \right) \right.
\\ & \ \ \ \ \ \ \ \ \ \ \ \left.- \frac{|2B|  }{|M_{24}|}\left( \chi_i(2B) - \frac{\chi_j(2B)}{\dim M^{(M_{24})}_j} \dim M^{(M_{24})}_i \right)\right)
\end{align*}
as $n \to \infty$, where $j = 1$ if $n$ is even and $j=2$ if $n$ is odd.
\end{corollary}

Finally, we show that our results apply to a sequence $V_n$ of $G$-modules for some vertex operator algebra $V=\bigoplus V_n$ and $G$ a finite group of automorphisms of $V$. More specifically, we have the following result.


\begin{theorem} Let $V = \bigoplus V_n$ be a holomorphic, $C_2$-cofinite, and self-dual vertex operator algebra. Let $G$ be a finite group of automorphisms of $V$. Let $g\in G$ and denote by $V(g)$ the unique (up to equivalence) $g$-twisted sector of $V$ . If the conformal weights satisfy $\rho(V(g))>\rho(V)$ for all $g\neq e$, then the sequence $(V_n)$ of $G$-modules has dominant identity trace. Consequently, $V_n$ tends to a multiple of the regular representation as $n \to \infty$.
\end{theorem}

We note here that the assumption on the conformal weights is conjectured by \cite{moller} to always hold for $V$ a holomorphic, $C_2$-cofinite vertex operator algebra, and $G$ a finite group of automorphism of $V$. Section $3.2$ provides cases where this assumption are known to hold.

\begin{acknowledgement} {We are thankful to John Duncan and Ken Ono for suggesting the problem and for their indispensable advice and guidance throughout the process. We are also grateful to Scott Carnahan, Robert Griess, and Hannah Larson for helpful conversations, and to an anonymous referee for comments that substantially improved the paper.} \end{acknowledgement}
\section{Proofs of Theorems 1.1, 1.2, and 1.4}
\subsection*{Proof of Theorem 1.1}
Let $\text{m}_i(n)$ be the multiplicity of $M_i$ in $K_n$ so that
$$K_n = \text{m}_1(n)M_1\oplus \text{m}_2(n)M_2\oplus\cdots\oplus \text{m}_s(n)M_s.$$
Let $\chi_i$ be the irreducible character of $M_i$. By the usual orthogonality of characters,
\begin{equation} \label{min}
\text{m}_i(n)=\dfrac{1}{|G|}\sum_{g\in G}\overline{\chi_i(g)}c_g(n).
\end{equation}
Since $(K_n)$ has dominant identity trace, we see that the formula for $\text{m}_i(n)$ is dominated by the term corresponding to $g=e$, and we recover the asymptotic given in Theorem \ref{thm1}.
\hfill$\square$

\subsection*{Proof of Theorems 1.2 and 1.4}

Let $\mathcal{C}$ be the set of conjugacy classes of $G$. We define a partial ordering on the elements of $\mathcal{C}$ as follows: $[h] < [g]$ if and only if $c_h(n) = o(c_g(n))$ as $n \to \infty$, and $[g] = [h]$ if and only if $c_g(n) \sim kc_h(n)$ for some constant $k$ as $n\to \infty$. In other words, we order the elements of $\mathcal{C}$ by increasing order of growth of the corresponding $c_g(n)$'s.

Let $L_1$ be the regular representation of $G$ and decompose $K_{n}$ into the direct sum of representations
\begin{equation}	\label{Kn} K_{n} = r_1(n)L_1 \oplus K_{n}^{(1)}
\end{equation}
where $r_1(n)$ is a nonnegative integer which is as large as possible. (Thus $K_{n}^{(1)}$ is the non-free part of $K_{n}$).  Then by definition,

\begin{align}
r_1(n) & := \min_{1\leq j\leq s}\left\{ \left\lfloor \dfrac{\text{m}_j(n)}{\dim M_j} \right\rfloor \right\} \nonumber\\
  &= \min_{1\leq j\leq s} \left\{ \left\lfloor \dfrac{1}{|G|}\left( c_{e}(n) + \sum_{[g] \in \mathcal{C}, g \neq e} \dfrac{|[g]|\overline{\chi_j(g)}c_{g}(n)}{\dim M_j} \right) \right\rfloor \right\}
\end{align}

Let
\begin{equation} \label{C2}
\mathcal{C}_2 = \{ [g] \in \mathcal{C} : [h] \leq [g] < [e] \text{ for all } h\in G\setminus \{e\} \}.
\end{equation} That is, $\mathcal{C}_2$ is the collection of conjugacy classes in $G$ whose corresponding $c_g(n)$'s has the second fastest growth rate (with $c_e(n)$ having the fastest growth by assumption).
As $n \to \infty$, the dominant terms in the sum in (2.3) are the terms corresponding to the conjugacy classes in $\mathcal{C}_2$. Therefore when $n$ is sufficiently large, we may find $r_1(n)$ by finding a $j$ which minimizes
	\[ \dfrac{\sum\limits_{[g] \in \mathcal{C}_2} |[g]|\overline{\chi_j(g)}c_g(n)}{\dim M_j}, \]
or equivalently, a $j$ which minimizes
\[ \dfrac{\sum\limits_{[g] \in \mathcal{C}_2} |[g]|\overline{\chi_j(g)}\text{sgn}(c_g(n))}{\dim M_j}. \] 
Let $j_1$ be one such $j$, which exists since the signs are independent of $n$, so that
	\[ j_1 \in J_1 := \left\{ j \in \{1,\ldots,s\} : j \text{ minimizes } \left\lfloor \dfrac{\text{m}_j(n)}{\dim M_j}\right\rfloor \text{ for } n \text{ sufficiently large}  \right\}. \]
	
Thus, if we write $K_{n}^{(1)}$ in terms of irreducible representations of $G$, say
	\[ K_{n}^{(1)} = \text{m}_1^{(1)}(n)M_1 \oplus \cdots \oplus \text{m}_s^{(1)}(n)M_s, \]
then the multiplicity functions for the non-free part have the following exact formula:
\begin{equation} \label{mi1n}
\text{m}_i^{(1)}(n) = \text{m}_i(n) - \left\lfloor\dfrac{\text{m}_{j_1}(n)}{\dim M_{j_1}}\right\rfloor \dim M_i.
\end{equation}
Using (\ref{min}), we have 
	\[ \text{m}_i^{(1)}(n) \sim \dfrac{1}{|G|} \sum_{[g] \in \mathcal{C}} |[g]|f_i^{(1)}(g)c_g(n) \]
as $n \to \infty$, where
	\[ f_i^{(1)}(g) := \left(\overline{\chi_i(g)} -  \dfrac{\dim M_i}{\dim M_{j_1}}\overline{\chi_{j_1}(g)}\right). \]
Note that $f_i^{(1)}(e) = 0$ so that the dominant terms in the asymptotic formula for $\text{m}_i^{(1)}(n)$ are the terms corresponding to the conjugacy classes in $\mathcal{C}_2$. More precisely,
	\[ \text{m}_i^{(1)}(n) \sim \dfrac{1}{|G|} \sum_{[g] \in \mathcal{C}_2} |[g]| f_i^{(1)}(g)c_{g}(n) \]
and this proves Theorem \ref{thm3}.

Note that when $i \in J_1$, as $n \to \infty$ we have
	\[ \left\lfloor\dfrac{\text{m}_{j_1}(n)}{\dim M_{j_1}}\right\rfloor = \left\lfloor\dfrac{\text{m}_{i}(n)}{\dim M_{i}}\right\rfloor \]
so that the multiplicity function \ref{mi1n} is equal to
	\[ \text{m}_i^{(1)}(n) = \text{m}_i(n) - \left\lfloor\dfrac{\text{m}_{i}(n)}{\dim M_{i}}\right\rfloor \dim M_i \leq \dim M_i. \]
Hence $K_{n}^{(1)}$ tends to a multiple of a representation $L_2$ whose irreducible components do not contain $M_i$ for $i \in J_1$.

Similar to (\ref{Kn}), we can then write
	\[ K_{n}^{(1)} = r_2(n)L_2 \oplus K_{n}^{(2)} \]
where $r_2(n)$ is a nonnegative integer that is as large as possible. We can repeat the arguments as before to show that there exist a set $J_3$ and a representation $L_3$ whose irreducible components do not contain $M_i$ for $i \in J_1 \cup J_2$, such that $K_n^{(2)}$ tends to a multiple of $L_3$.
Proceeding inductively, this gives us the decomposition
	\[ K_{n} = r_1(n)L_1 \oplus r_2(n)L_2 \oplus \cdots \oplus r_t(n)L_t \oplus L_{\epsilon}(n), \]
where the $L_j$ are expressed in terms of fewer and fewer irreducible representations of $G$, and where $L_{\epsilon}(n)$ is a representation of $G$ with bounded multiplicity functions. This proves Theorem \ref{thm2}.
\hfill $\square$

\section{Applications}
\subsection{Umbral moonshine}

The umbral moonshine conjecture, proven in \cite{proof}, states that for a Niemeier root system $X$ and setting $m=m^X$ where $m^X$ is the Coxeter number of any simple component of $X$, there is a naturally defined bi-graded infinite-dimensional representation of $G^X$ $$K^X=\bigoplus\limits_{r\in I^X} \bigoplus\limits_{\substack{D\in \mathbb{Z}, D\leq0 \\ D=r^2\pmod{m} }} K^X_{r, -D/4m}$$ such that the graded trace functions  $$H_{g,r}^X(\tau)=-2q^{-1/4m}\delta_{r,1}+\sum\limits_{\substack{D\in \mathbb{Z}, D\leq0 \\ D=r^2\pmod{m} }} tr(g|K^X_{r, -D/4m})q^{-D/4m}$$ for $r\in I^X$ are components of \textit{vector-valued} mock modular forms $H_g^X$. Here $I^X$ depends on the types of components in the root system \cite{proof}. 

Convenient expressions for mock modular forms can be found using Rademacher sums which are essentially regularized Poincar\'e series. The theory of Rademacher sums, dating back to the 1930s and originally establishing a conditionally convergent expression for the normalized modular $j$-invariant, has been generalized to apply to modular and mock modular forms of various weights and various subgroups of $SL_2(\mathbb{R})$ \cite{Hans, MR0344196, DF, Rademacher}. Cheng and Duncan considered Niebur's method for constructing Rademacher sums in weight $1/2$ and verified that it produces the functions appearing in Mathieu moonshine \cite{MR3021323}.  
As an example, we examine the case of Mathieu moonshine. Following Cheng and Duncan \cite{Rademacher}, let $n_g$ be the order of $g$, let $\epsilon$ be the multiplier system for the Dedekind $\eta$ function, and let $\rho_g$ be a character specified by the minimal cycle length in the cycle shape of $g$.
Then the mock modular form $H_g^{(M_{24})}$ is defined as:
\begin{align*}H_g^{(M_{24})}(\tau)&=-2R^{[-1/8]}_{\Gamma_0(n_g), \rho_g\epsilon^{-3}, 1/2}(\tau)\\ &=-2q^{-1/8}+2\sum_{n>0} c_{\Gamma_0(n_g), \rho_g\epsilon^{-3}, \frac{1}{2}} q^{n-1/8}.\end{align*}
The coefficients can be expressed in terms of Kloosterman sums weighted by Bessel functions which have known asymptotics for which our theory applies.

 \subsection*{Proof of Corollary 1.5} From the Rademacher sum formulation, we find that the $n$th coefficient of $H_g^{(M_{24})}$ is
	\[ c_g(n) =-4\pi\sum\limits_{\substack{c>0 \\ 0\leq d <c \\ (c,d)=1}}e\left(n\frac{d}{c}-\frac{3s(d,c)}{2}\right)e\left(-\frac{cd}{n_gh_g}\right) \frac{1}{c(8n-1)^{\frac{1}{4}}}I_{\frac{1}{2}}\left(\frac{\pi}{2c}(8n-1)^{\frac{1}{2}}\right)	.\]
From the asymptotics of the $I$-Bessel function \[ I_\frac{1}{2}(x) \sim \frac{e^x}{\sqrt{2\pi x}}, \] and upon isolating the dominant term of this sum, we get
	\[ c_g(n) = \text{sgn}(c_g(n))\frac{4}{n_g^{1/2}\sqrt{8n-1}} \text{exp}\left(\frac{\pi\sqrt{8n-1}}{2n_g}\right) +o \left(\text{exp}\left(\frac{\pi\sqrt{8n-1}}{2n_g}\right)\right) .\]
Thus $\mathcal{C}_2 = \{2A, 2B\}$. Also, from the Rademacher sum formulation, we find that $\text{sgn}(c_{2A}(n)) = (-1)^n $ and $\text{sgn}(c_{2B}(n)) = (-1)^{n+1}$. Thus, from the character table of $M_{24}$, the $j$ that minimizes 
	\[ \dfrac{\sum\limits_{g \in \mathcal{C}_2} |[g]|\chi_j(g)\text{sgn}(c_g(n))}{\dim M_j^{(M_{24})}} \]
is $j =1 $ when $n$ is even, and $j =2 $ when $n$ is odd. (Note that $\chi_j(2A)$ and $\chi_j(2B)$ are real-valued for any $j$.) Thus Theorem \ref{thm3} in this case becomes:
\begin{align*}
	\text{m}_i'(n) &\sim (-1)^{n+1}\frac{2\sqrt{2} \, e^{\frac{\pi}{4}\sqrt{8n-1}}}{\sqrt{8n-1}}\left(\frac{|2A|  }{|M_{24}|}\left( \chi_i(2A) - \frac{\chi_j(2A)}{\dim M^{(M_{24})}_j}\dim M^{(M_{24})}_i \right) \right.
\\ & \ \ \ \ \ \ \ \ \ \ \ \left.- \frac{|2B|  }{|M_{24}|}\left( \chi_i(2B) - \frac{\chi_j(2B)}{\dim M^{(M_{24})}_j} \dim M^{(M_{24})}_i \right)\right). 
\end{align*}\hfill $\square$

Similar to the Mathieu moonshine case, the vector-valued mock modular forms appearing in umbral moonshine, denoted $H_g^{(\ell)}$ for lambency $\ell$, can also be conveniently expressed in terms of \textit{ vector-valued} Rademacher sums according to \[H_g^{(\ell)}=-2R_{\Gamma_0(n_g),\psi^{(\ell)}\rho_g^{(\ell)}, 1/2}\] which have similar asymptotics to the usual Rademacher sums \cite{Rademacher, proof}.

\subsection{Vertex Operator Algebras}
In this section, we consider vertex operator algebras. We refer to \cite{FrenkelBenzvi} or \cite{MR996026} for basic definitions in the theory of vertex operator algebras. Let $V = \bigoplus V_n$ be a vertex operator algebra and let $G$ be any finite group of automorphisms of $V$. We show that if $V$ and $G$ satisfy certain natural conditions (to be enumerated shortly), then the sequence $(V_n)$ of $G$-modules has dominant identity trace. Consequently, Theorem 1.1 applies to such a sequence of modules.

The natural conditions for $V$ that we will assume are the following: holomorphic, $C_2$-cofinite, and self-dual. By Theorem 2 of \cite{donglimason}, if $g$ is an automorphism of $V$ of finite order, then $V$ possesses a unique $g$-twisted sector (up to equivalence), which is denoted $V(g)$. We shall assume another condition for $V$ and $G$ concerning the conformal weights of the twisted sectors; we assume that the conformal weight $\rho(V)$ of the untwisted sector $V$ is strictly less than the conformal weights $\rho(V(g))$ of the other twisted sectors $V(g)$ for $g \in G\setminus \{e\}$. This minimality condition is conjectured to always hold when $V$ is holomorphic and $C_2$-cofinite, and when $G$ is a finite group of automorphisms of $V$ (cf. \ Conjectures 1.1 and 2.2 in \cite{moller}). We have the following proposition.

\vspace{0.1in}

\begin{thmV} Let $V = \bigoplus V_n$ be a holomorphic, $C_2$-cofinite, and self-dual vertex operator algebra. Let $G$ be a finite group of automorphisms of $V$. Let $g\in G$ and denote by $V(g)$ the unique (up to equivalence) $g$-twisted sector of $V$. If the conformal weights satisfy $\rho(V(g))>\rho(V)$ for all $g\neq e$, then the sequence $(V_n)$ of $G$-modules has dominant identity trace. Consequently, $V_n$ tends to a multiple of the regular representation as $n \to \infty$.
\end{thmV}

\vspace{0.1in}

There are examples for which the assumptions in Theorem 1.6 are known to hold. For instance, if $W$ is a holomorphic and $C_2$-cofinite vertex operator algebra, and if $k \in \mathbb{Z}$, then our assumptions hold for $V = W^{\otimes k}$ and $G \leq S_k \leq \text{Aut}(V)$ (cf.\ Proposition $4.1$ of \cite{moller}). Other examples come from lattice vertex operator algebras. Given an even, unimodular, and positive-definite lattice $L$, our assumptions also hold for $V = V_L$ (the lattice vertex operator algebra associated to $L$) and 
 $G\leq \text{Aut}(V_L)$ (cf.\ Proposition $4.2$ of \cite{moller}). We refer the reader to Section $4.2$ of \textit{loc.~ cit.} for more on the structure of $\text{Aut}(V_L)$. In particular, for these examples of $V$ and $G$, the sequence $(V_n)$ of $G$-modules has dominant identity trace, and thus, Theorem $1.1$ applies.

\vspace{0.1in}

\subsection*{Proof of Theorem 1.6} Let $V := V(e)$, and write $V(g)$ as follows:
\[ V(g)=\bigoplus_{n=0}^{\infty} V(g)_{\frac{n}{\text{ord}(g)}+\rho(V(g))}.\]
If $h \in G$ commutes with $g$, then $h$ induces an action on $V(g)$ (which is well-defined up to a scalar factor). If $c$ is the central charge of $V$, then the twisted trace functions are the following power series:
\begin{equation}\label{eqa}Z(g,h;\tau):= \sum_{n=0}^\infty \text{tr}(h\mid V(g)_{\frac{n}{\text{ord}(g)}+\rho(V(g))})q^{\frac{n}{\text{ord}(g)}+\rho(V(g))-\frac{c}{24}}.\end{equation}
In particular, for all $h \in G$, the coefficients of $Z(e,h;\tau)$ encode the traces of $h$ on the homogenous subspaces of $V$. We denote by $c_h(n)$ these coefficients, i.e.,
\[ Z(e,h;\tau) = \sum_{n=0}^{\infty} c_h(n) q^{n+\rho(V) - \frac{c}{24}}. \]
To prove the proposition, we need to show that $c_h(n) = o(c_e(n))$ for $h\neq e$ as $n \rightarrow \infty$. Now, since $V$ is holomorphic and $C_2$-cofinite, we know from \cite{donglimason} that $Z(g,h;\tau)$ is holomorphic in $\mathbb{H}$, and moreover, if $\gamma =\left(\begin{smallmatrix} a & b \\
c & d \end{smallmatrix}\right) \in SL_2(\mathbb{Z})$, then
\begin{equation}\label{eqb}Z(g,h;\gamma\tau)= \sigma(g,h,\gamma) Z(g^ah^c,g^bh^d;\tau)\end{equation}
for some constant $\sigma(g,h,\gamma)$. Furthermore, because $V$ is also self-dual, the invariance subgroups of these graded trace functions are congruence subgroups \cite{DongLi}. 

Let $h\neq e$. Choose $\gamma =\left(\begin{smallmatrix} a & b \\
c & d \end{smallmatrix}\right) \in SL_2(\mathbb{Z})$ such that $c$ is not a multiple of the order of $h$. By (\ref{eqb}), the expansion of $Z(e,h;\tau)$ at the cusp $\gamma\infty$ is $\sigma(e,h,\gamma) Z(h^c, h^d;\tau)$. Since $h^c \neq e$, we have $\rho(V(h^c))>\rho(V)$ by our assumption on conformal weights. So by (\ref{eqa}), the order of the pole of $Z(e,e;\tau)$ at $\tau = \gamma\infty$ is strictly greater than the order of the pole of $Z(e,h;\tau)$ at $\tau = \gamma\infty$. At any other cusp $s$, the order of the pole of $Z(e,e;\tau)$ at $\tau = s$ is greater than or equal to the order of the pole of $Z(e,h;\tau)$ at $\tau = s$. 

We will use this comparison of the orders of poles at the cusps together with the fact that these twisted trace functions are modular functions on congruence subgroups to compare the asymptotic growths of the coefficients. That is, since the trace functions are modular they can be expressed in terms of Rademacher sums and their coefficients can be expressed as Rademacher series \cite{DF}. In particular, from the asymptotic expressions in Section $4.2$ of \cite{DF} we can read the growth of $c_h(n)$ at each of the different cusps. This allows us to see how the growth of the coefficients depends on the order of the pole, and we find that the Bessel function (cf. $(4.2.1)$ of \cite{DF}) is largest when the order of the pole is largest. Thus we have that $c_h(n) = o(c_e(n))$ and so the sequence $(V_n)$ of $G$-modules has dominant identity trace and Theorem \ref{thm1} applies. 
\hfill$\square$

\bibliography{MathieuBibtex}{}
\bibliographystyle{plain}

\end{document}